\documentclass[11pt,a4paper]{article}
\usepackage{caption2}%---change the caption of figures
\usepackage{CJK}
\usepackage{tabls}
\usepackage{bm}
\usepackage{multirow}
\usepackage{amssymb}
\usepackage{amsfonts}
\usepackage{mathrsfs}
\usepackage{empheq}
\usepackage{amsmath}
\usepackage{amsthm}
\usepackage{graphicx}
\usepackage{color}
\usepackage{indentfirst}
\usepackage{hyperref}%pour PDFLaTex
\usepackage{float}
\usepackage{cases}
\usepackage[titletoc]{appendix}
\usepackage{bm}

\usepackage[colorinlistoftodos,english]{todonotes}
\allowdisplaybreaks%

\def\i1n{i=1,\cdots,n}
\def\j1n{j=1,\cdots,n}
\def\ij1n{i,j=1,\cdots,n}

\def \i{\mathrm i}

 \numberwithin{equation}{section}
\theoremstyle{definition}

 \newtheorem{thm}{\indent Theorem}[section]

 \newtheorem{rem}{\indent Remark}[section]

\theoremstyle{definition}
\theoremstyle{theorem}
\theoremstyle{lemma}

\newcommand{\be}{\begin{equation}}
\newcommand{\ee}{\end{equation}}
\newcommand{\beq}{\begin{equation*}}
\newcommand{\eeq}{\end{equation*}}

\begin{document}

\begin{CJK*}{GB}{gbsn}

\title{Generalized persistence of entropy weak solutions for system of hyperbolic conservation laws}
\author{Yi Zhou\thanks{School of Mathematical sciences, Fudan
        University, Shanghai 200433, P. R. China, {\it email: yizhou@fudan.edu.cn.}}}
\date{\today}
\maketitle

\begin{abstract}
Let $u(t,x)$ be the solution to the Cauchy problem of a scalar conservation law in one space dimension. It is well known that even for smooth initial data the solution can become discontinuous in finite time and global entropy weak solution can best lie in the space of bounded total variations. It is impossible that the solutions belong to ,for example ,$H^1$ because by Sobolev embedding theorem $H^1$ functions are H$\mathrm{\ddot{o}}$lder continuous. However, we note that from any point $(t,x)$ we can draw a generalized characteristic downward which meets the initial axis at $y=\alpha (t,x)$. if we regard $u$ as a function of $(t,y)$, it indeed belongs to $H^1$ as a function of $y$ if the initial data belongs to $H^1$. We may call this generalized persistence (of high regularity) of the entropy weak solutions. The main purpose of this paper is to prove some kinds of generalized persistence (of high regularity) for the scalar and $2\times 2$ Temple system of hyperbolic conservation laws in one space dimension .
\end{abstract}
\textbf{Keyword:}  quasilinear hyperbolic system, Cauchy problem, entropy weak solution. vanishing viscosity method.

\section{Introduction}%----------------------------------------------------------------introduction
 The Cauchy problem for system of conservation laws in one space
 dimension takes the form
 \begin{equation}\label{1.1}%--------------------------------------------------------------1.1
 u_t + f(u)_x = 0,
 \end{equation}
 \begin{equation}\label{1.2}
 u(0,x) = u_0(x).
 \end{equation}
 Here $u=(u_1,\cdots ,u_n)$ is the vector of conserved quantities,
 while the components of $f=(f_1,\cdots ,f_n)$ are the fluxes. We
 assume that the flux function $f:R^n\to R^n$ is smooth and that
 the system is  hyperbolic,i.e. at each point $u$ the
 Jacobian matrix $A(u)=\nabla f(u)$ has $n$ real eigenvalues
 \begin{equation}\label{1.3}
 \lambda_1(u) , \cdots  ,\lambda_n(u)£¬
 \end{equation}
 and a bases of right and left eigenvectors
  $r_i(u),l_i(u)$, normalized so that
 \begin{equation}\label{1.4}%------------------------------------------------------------1.4
 l_i\cdot r_j = \delta_{ij},
 \end{equation}
 where $\delta_{ij}$ stands for  Kroneker's symbol. We make an assumption that all the eigenvalues and eigenvectors are smooth functions of $u$, which in particular hold when the eigenvalues are all distinct, i.e. the system is strictly hyperbolic.

 It is well known that the solution can develop singularities in
 finite time even with smooth initial data, see Lax \cite{Lax2},John \cite{John}, Li \cite{Li} . Therefore, global solutions can only be constructed within a
 space of discontinuous functions. Global weak solutions to the
 Cauchy problem  is a subject of a large literature,
 notably, Lax \cite{Lax1}, Glimm \cite{Glimm}, DiPerna \cite{DiPerna}, Dafermos \cite{Dafermos1}, Bressan \cite{Bressan}, Bianchini and Bressan \cite{Bianchini}. We
 refer to the classical monograph  of Dafermos \cite{Dafermos2} for references.

  In the classical paper of  Kruzkov \cite{Kruzkov}, global entropy weak
 solutions to a scalar equation are constructed by a vanishing
 viscosity method. That is , the entropy weak solutions of the
 hyperbolic equation  actually coincide with the limits of
 solutions to the parabolic equation
 \begin{equation}\label{1.5}%----------------------------------------------------------------1.5
 u_t + f(u)_x = \varepsilon u_{xx}
\end{equation}
 by letting the viscosity coefficients $\varepsilon\to 0$. The same
  result is also proved for $n\times n$ strictly hyperbolic systems in a cerebrated paper of Bianchini and Bressan \cite{Bianchini} for small BV initial data.

Although the one dimension theory of systems of hyperbolic conservation laws has by now quite matured,  the multi-dimensional  problem is still very challenging  except
for the scalar case. In recent years, much progress has been made to understand the formation of shocks for the compressible Euler equations for small initial data, see Sideris \cite{Sideris} , Christodoulou \cite{Christodoulou}. However, the problem of constructing entropy weak solutions beyond the time of shock formation is still largely open and even so in the radial symmetric case. The main difficulty is that on the one hand the weak solution can best lie in BV, on the other hand, the BV space is not a scaling invariant space for the system. Especially, in n space dimensions, $\dot{W}^{1,n}$ is the critical space. This motivates us to study systems of hyperbolic conservation laws in one space dimension for $W^{1,p},(1<p<+\infty) $ data as a first step towards the multidimensional problem.

Let $u(t,x)$ be the solution to the Cauchy problem of a scalar conservation law in one space dimension. It is well known that even for smooth initial data the solution can become discontinuous in finite time and  global entropy weak solution can best lie in the space of bounded total variations. It is impossible that the solution belong to $W^{1,p}, (1<p<+\infty)$ because by Sobolev embedding theorem $W^{1,p}$ functions are Holder continuous. However, we note that from any point $(t,x)$ we can draw a generalized characteristic downward which meets the initial axis at $y=\alpha (t,x)$. if we regard u as a function of $(t,y)$, it indeed belongs to $W^{1,p}$ as a function of y if the initial data belongs to $W^{1,p}$. We may call this generalized persistence (of high regularity) of the entropy weak solutions. The main purpose of this paper is to prove some kind of generalized persistence of $W^{1,p}$  regularity of entropy weak solutions for  $2\times 2$ Temple system of  hyperbolic conservation laws in one space dimension. Some interesting hyperbolic system arising in applications which satisfies the Temple condition can be found in Serre \cite{Serre}.

 Our main theorem can be stated as follows:
\begin{thm}\label{thm1.1}
Consider the Cauchy problem \eqref{1.1}, \eqref{1.2} for systems of two conservation laws. Suppose that
the $2\times 2$ matrix $A(u)$ is hyperbolic, smoothly depending on u and possess a complete sets of smooth eigenvalues and eigenvectors  as well as  two global Riemann invariants.
 Suppose that the Temple condition
\begin{equation}\label{1.6}
<l_i, r_{ju}r_j>\equiv 0, \forall i\not=j.
\end{equation}
is satisfied. Suppose furthermore that
\begin{equation}
|u_0|_{L^\infty(R)}\le D<+\infty
\end{equation}
and there exists $1<p\le +\infty$ such that
\begin{equation}\label{1.8}
|u_0'|_{L^p(R)}=M<+\infty .
\end{equation}
then, the Cauch problem \eqref{1.1},\eqref{1.2} admits a global entropy weak solution which can be represented as
\begin{equation}
u=U( W_1(t, \alpha_1(t,x)), W_2(t, \alpha_2(t,x)) ),
\end{equation}
where $U$ is a smooth function of Riemann invariants $W_1, W_2$ ,$\alpha_1(t,x), \alpha_2(t,x)$ are locally bounded monotone increasing function of $x$  and $W_1(t,\alpha), W_2(t,\alpha)$ are Holder continuous function of $\alpha$ and moreover
\begin{equation}
|\partial_{\alpha}W_1(t,\alpha )|_{L^p(R)}+|\partial_{\alpha}W_2(t,\alpha )|_{L^p(R)}\lesssim M .
\end{equation}
\end{thm}
 Theorem1.1 will be proved by a vanishing viscosity approach.
 \begin{rem}Theorem 1.1 is also true for the initial boundary value problems with periodic boundary conditions, with \eqref{1.8} replaced by
 \begin{equation}
 |u_0'|_{L^p(T)}=M<+\infty .
 \end{equation}
  The same proof applies.
 \end{rem}

\begin{rem}With additional assumption \eqref{1.8}, Theorem 1.1 gives an alternative proof of global existence of entropy weak solution for $2\times 2$ Temple system without using the so called compensated compactness method.
\end{rem}

 This paper is organized as follows:
 In section2, we will discuss generalized persistence of a scalar conservation law in one space dimension in various high regularity spaces. In section 3, we will discuss related problem for a scalar conservation law in multi-dimensions. Finally, in section 4, we will discuss the $2\times 2$ Temple system in one space dimension and prove our main result.

 Notations: Let $f(x)$ be a scalar or vector function of $x\in R$, we denote
\begin{equation}
|f|_{\dot{W}^{1,p}(R)}=|f'|_{L^p(R)}
\end{equation}
\begin{equation}
|f|_{W^{1,p}(R)}=|f|_{\dot{W}^{1,p}(R)}+|f|_{L^p(R)}
\end{equation}
and $\dot{H}^1=\dot{W}^{1,2}$.
We denote $A\lesssim B$, if there exist a positive constant C such that $A\le CB$.

\section{Scalar equation in one space dimension}
We consider the following Cauchy problem for a scalar conservation law in one space dimension:
\begin{equation}\label{2.1}
u(t,x)_t+f(u(t,x))_x=0,
\end{equation}
\begin{equation}\label{2.2}
u(0,x)=u_0(x),
\end{equation}
where $u_0$ is a suitably smooth function.
It is well known that the global solution is  the limit of the viscous approximations
\begin{equation}\label{2.3}
u^\varepsilon(t,x)_t+f(u^\varepsilon(t,x))_x=\varepsilon u^{\varepsilon}_{xx},
\end{equation}
\begin{equation}\label{2.4}
u^\varepsilon(0,x)=u_0(x).
\end{equation}
By maximum principle, we have
\begin{equation}\label{2.6}
|u^\varepsilon|_{L^\infty (R^+\times R)}\le |u_0|_{L^\infty }.
\end{equation}
We write
\begin{equation}\label{2.5}
u^\varepsilon(t,x)=U^\varepsilon(t,\alpha^\varepsilon (t,x)),
\end{equation}
for simplicity of notation, here we denote $U^\varepsilon(t,\alpha^\varepsilon (t,x))$ just by $U(t,\alpha (t,x))$ .
Plug \eqref{2.5} to \eqref{2.3}, we get
\begin{equation}
U(t,\alpha )_t-\varepsilon \alpha_x^2U(t,\alpha )_{\alpha \alpha}=-U(t,\alpha )_\alpha(\alpha_t+f'(u)\alpha_x-\varepsilon \alpha_{xx}).
\end{equation}
We take $\alpha^\varepsilon(t,x)$ to be the solution to the following Cauchy problem
\begin{equation}\label{2.7}
\alpha^\varepsilon_t+f'(u^\varepsilon )\alpha^\varepsilon_x-\varepsilon \alpha^\varepsilon_{xx}=0,
\end{equation}
\begin{equation}\label{2.8}
\alpha^\varepsilon (0,x)=x,
\end{equation}
then, $U^\varepsilon (t,\alpha )$ will satisfy
\begin{equation}\label{2.9}
U^\varepsilon(t,\alpha )_t-\varepsilon \Theta^{\varepsilon 2}U^\varepsilon(t,\alpha )_{\alpha \alpha}=0,
\end{equation}
\begin{equation}\label{2.10}
U^\varepsilon (0,\alpha )=u_0(\alpha ),
\end{equation}
where we denote $\Theta^\varepsilon =\alpha^\varepsilon_x$. By \eqref{2.7}, \eqref{2.8}, we get
\begin{equation}\label{2.11}
\Theta^\varepsilon_t+(f'(u^\varepsilon )\Theta^\varepsilon)_x-\varepsilon \Theta^\varepsilon_{xx}=0,
\end{equation}
\begin{equation}\label{2.12}
\Theta^\varepsilon (0,x)=1.
\end{equation}
Then by maximum principle, $\Theta^\varepsilon$ is a positive function, and moreover, by \eqref{2.7},\eqref{2.8}
\begin{equation}\label{2.13}
x+M_1t\le \alpha^\varepsilon (t,x)\le x+M_2t
\end{equation}
where $M_1=\inf_{|u|\le |u_0|_{L^{\infty}}}(-f'(u))$, $M_2=\sup_{|u|\le |u_0|_{L^{\infty}}}(-f'(u)) $
This is because $x+M_1t$   and $x+M_2t$ is respectively a subsolution and supsolution of \eqref{2.7} ,\eqref{2.8}.
Now, form \eqref{2.9}, it is easy to see $V^\varepsilon(t,\alpha)=U^\varepsilon(t,\alpha)_\alpha$ satisfies
\begin{equation}\label{2.14}
V^\varepsilon(t,\alpha )_t-\varepsilon \big( \Theta^{\varepsilon 2}V^\varepsilon(t,\alpha )_{\alpha }\big)_\alpha=0.
\end{equation}
Then, it is easy to get the following series of estimates
\begin{equation}
|U^\varepsilon(t,\cdot )_{\alpha\alpha}|_{L^1(R)}\le |u_0''|_{L^1(R)}
\end{equation}
and for any $1\le p \le \infty$
\begin{equation}
|U^\varepsilon(t,\cdot )_{\alpha}|_{L^p(R)}\le |u_0'|_{L^p(R)}.
\end{equation}
Upon taking a subsequence, $U^\varepsilon (t,\alpha)$ converges to $ U(t,\alpha)$ and $\alpha^\varepsilon (t,x)$ converges to $\alpha (t,x)$.
Then $u(t,x)=U(t, \alpha (t,x))$ is the solution to the Cauchy problem \eqref{2.1},\eqref{2.2}, We see immediately that
$U(t,\alpha)_\alpha$ is a function of bounded total variation for the variable $\alpha$ provided that  $u_0'$ is a function of bounded total variation  and
$U(t,\alpha)_\alpha$ is an $L^p$ ($1<p\le \infty $) function for the variable $\alpha$  provided that  $u_0'$ is an $L^p$ function. Thus,
$U(t,\alpha)$ can be  much smoother then $u(t,x)$.
We summarize our result in the following theorem
\begin{thm}\label{thm2.1}
Let
\begin{equation}
u_0\in L^\infty ,
\end{equation}
then the global entropy solution to system \eqref{2.1} \eqref{2.2} can be represented as
\begin{equation}
u(t,x)=U(t,\alpha (t,x))
\end{equation}
where $\alpha (t,x)$ is a locally bounded monotone increasing function of x reprensating the generalized characteristics and $U(t,\alpha )$ as a function of $\alpha$ satisfies
\begin{equation}
|U(t,\cdot )_{\alpha}|_{BV(R)}\le |u_0'|_{BV(R)}
\end{equation}
and for any $1< p \le \infty$
\begin{equation}
|U(t,\cdot )_{\alpha}|_{L^p(R)}\le |u_0'|_{L^p(R)},
\end{equation}
provided that left hand side of the inequality is finite, i.e. $u_0$ is suitable smooth. In particular, $U(t,\alpha )$ is H$\mathrm{\ddot{o}}$lder continuous.
\end{thm}

\section{Scalar conservation law in multi-dimensions }
In this section, we consider the initial boundary value problem with periodic boundary conditions of a scalar conservation law in multi-dimensions
\begin{equation}\label{3.1}
u_t+\sum_{i=1}^n(f_i(u))_{x_i}=0,
\end{equation}
\begin{equation}\label{3.2}
u(0,x)=u_0(x), x\in T^n=[0,1]^n.
\end{equation}
As always, u is the limit of its viscous approximations:
\begin{equation}\label{3.3}
u^\varepsilon_t+\sum_{i=1}^n(f_i(u^\varepsilon ))_{x_i}=\varepsilon \Delta u^\varepsilon ,
\end{equation}
\begin{equation}\label{3.4}
u^\varepsilon (0,x)=u_0(x), x\in T^n,
\end{equation}
where $\Delta$ is the Laplacian operator in $T^n$.

Due to the multi-dimensional nature of the problem, there no longer exist a transformation $y=\alpha (t,x)$ like that in one space dimensions.
Therefore, in this section, we are limited to discuss regularity properties of solutions of the viscous approximations.

Let $\Theta^\varepsilon (t,x)$ be the solution to the initial boundary value problem with periodic boundary conditions of the following equation
\begin{equation}\label{3.6}
\Theta^\varepsilon_t+\sum_{i=1}^n(f'_i(u^\varepsilon )\Theta^\varepsilon )_{x_i}=\varepsilon \Delta \Theta^\varepsilon ,
\end{equation}
\begin{equation}\label{3.7}
\Theta^\varepsilon (0,x)=1.
\end{equation}
By maximum principle, $\Theta^\varepsilon $ is a positive function, moreover, integrating \eqref{3.6} in x yields
\begin{equation}\label{3.9}
|\Theta^\varepsilon(t,\cdot )|_{L^1(T^n)}=1.
\end{equation}

Let
$$v^\varepsilon (t, x)=u^{\varepsilon}_{x_1}(t,x),$$ then differentiate the equation \eqref{3.1} with respect to $x_1$ yields
\begin{equation}\label{3.8}
v^{\varepsilon}_t+\sum_{i=1}^n(f'_i(u^\varepsilon )v^\varepsilon )_{x_i}=\varepsilon \Delta v^\varepsilon ,
\end{equation}
\begin{equation}
v^\varepsilon (0,x)=u_{0x_1}(x).
\end{equation}
A simple computation shows
\begin{equation}
(\frac{v^\varepsilon}{\Theta^\varepsilon})_t+\sum_{i=1}^nf'_i(u^\varepsilon)(\frac{v^\varepsilon }{\Theta^\varepsilon })_{x_i}=\varepsilon (\Theta^\varepsilon)^{-2}\nabla \big(\Theta^{\varepsilon 2}\nabla (\frac{v^\varepsilon}{\Theta^\varepsilon})\big ).
\end{equation}
Then by maximum principle, we get
\begin{equation}\label{3.10}
\frac{|u_{x_1}^\varepsilon (t,x)|}{\Theta^\varepsilon (t,x)}\le |u_{0x_1}|_{L^\infty(T^n)},
\end{equation}
uniformly for all $(t,x)$. In a same way, we have
\begin{equation}\label{3.11}
\frac{|u_{x_i}^\varepsilon (t,x)|}{\Theta^\varepsilon (t,x)}\le |u_{0x_i}|_{L^\infty(T^n)}, \quad i=1,\cdots ,n
\end{equation}
uniformly for all $(t,x)$.

Let $1<p<\infty$, a further computation yields
\begin{eqnarray}
&&\big(\frac{|v|^p}{\Theta ^{p-1}}\big)_t+\sum_{i=1}^n\big (f_i'(u)\frac{|v|^p}{\Theta ^{p-1}}\big)_{x_i}
+\varepsilon\frac{4(p-1) }{p}\Theta\big(\nabla \big(\frac{|v|^{p/2}}{\Theta ^{p/2}}\big)\big)^2\\\nonumber
&&=\varepsilon \Delta  \big(\frac{|v|^p}{\Theta^{p-1}}\big ),
\end{eqnarray}
where we denote $v^\varepsilon$, etc by $v$, etc for simplicity of notation. Similar equality hold for $u^\varepsilon_{x_i}$.
we integrate the above equality to yield
\begin{equation}\label{3.12}
\int_{T^n}\frac{|u^\varepsilon_{x_i}(t,x)|^p}{(\Theta^\varepsilon (t,x)) ^{p-1}}dx+\varepsilon\frac{4(p-1) }{p}\int_0^t\int_{T^n}\Theta^\varepsilon\big(\nabla \big(\frac{|u_{x_i}^\varepsilon|^{p/2}}{(\Theta^\varepsilon )^{p/2}}\big)\big)^2dxd\tau=\int_{T^n}|u_{0x_i}|^pdx .
\end{equation}
\eqref{3.11}, \eqref{3.12} gives a kind of  $L^P$, ($1<p\le \infty$) bound for the derivatives of the solution. We notice that by Holder's inequality and \eqref{3.9}, for any $1\le p_1<p$
\begin{equation}
\big(\int_{T^n}\frac{|u^\varepsilon_{x_i}(t,x)|^{p_1}}{(\Theta^\varepsilon (t,x)) ^{p_1-1}}dx\big)^{1/p_1}\le \big(\int_{T^n}\frac{|u^\varepsilon_{x_i}(t,x)|^p}{(\Theta^\varepsilon (t,x)) ^{p-1}}dx\big)^{1/p}.
\end{equation}
Thus, when $\varepsilon\to 0$ ,$\frac{|u^\varepsilon_{x_i}(t,x)|^p}{(\Theta^\varepsilon (t,x)) ^{p-1}}dx$ will converge weakly to some measures $\Omega _{ip}(t)$ such that
its total measure is in an increasing order of p, and it in someway measures the high regularity  part of the solution.

\section{$2\times 2$ Temple system in one space dimension}
We consider the viscous approximations
\begin{equation}\label{4.1}
u_t^\varepsilon+A(u^\varepsilon)u^\varepsilon_x=\varepsilon u^\varepsilon_{xx},
\end{equation}
with initial conditions
\begin{equation}\label{4.2}
u^\varepsilon(0,x)=u_{0\varepsilon}(x)=J_{\varepsilon}\ast u_0,
\end{equation}
where $J_\varepsilon$ is the Friedriches mollifier.
We assume that there exists two Riemann invariants $R^\varepsilon_1=R_1(u^\varepsilon )$, $R^\varepsilon_2=R_2(u^\varepsilon )$, and we assume that we can also write $u^\varepsilon=U(R^\varepsilon_1,R^\varepsilon_2)$.
In the following calculations, we just denote $u^\varepsilon$ by $u$ and $R^\varepsilon_i$ by $R_i$ when there is no confusion of notation.
Taking $l_i(u)=\nabla_u R_i(u)$, we get
\begin{equation}
u_x=\sum_{i=1}^2R_{ix}r_i(u).
\end{equation}
Thus, we get
\begin{equation}\label{4.6}
u_{xx}=\sum_{j=1}^2R_{jxx}r_j+\sum_{j,k=1}^2R_{jx}R_{kx}r_{ju}r_k.
\end{equation}
Thus, taking inner product of \eqref{4.1} with $l_i$ and noting the Temple condition, we get
\begin{equation}\label{4.7}
R_{it}+\lambda_iR_{ix}=\varepsilon R_{ixx}+\varepsilon R_{ix}\sum_{j=1}^2D_{ij}R_{jx}
\end{equation}
where $D_{ij}$ are smooth function of $u^\varepsilon$.
Thus,
\begin{equation}\label{4.7}
R_{it}+\tilde{\lambda_i}R_{ix}=\varepsilon R_{ixx}
\end{equation}
where
\begin{equation}
\tilde{\lambda_i}=\lambda_i-\varepsilon \sum_{j=1}^2 D_{ij}R_{jx}.
\end{equation}
Then, by maximum principle
\begin{equation}
|R_i|_{L^\infty}\le |R_i(u_{0\varepsilon})|_{L^\infty(R)}\lesssim 1.
\end{equation}
Thus, we have
\begin{equation}
|u^\varepsilon|_{L^\infty}\lesssim 1.
\end{equation}
We write
\begin{equation}\label{4.8}
R_i(t,x)=W_i(t, \alpha_i(t,x)),
\end{equation}
then \eqref{4.7} becomes
\begin{equation}
W_{it}-\varepsilon\alpha_{ix}^2W_{i\alpha \alpha}=-W_{i\alpha}\big(\alpha_{it}+\tilde{\lambda_i}\alpha_{ix}-\varepsilon\alpha_{ixx}\big).
\end{equation}
Thus, we get
\begin{equation}\label{4.10}
W^\varepsilon_{it}-\varepsilon(\alpha^\varepsilon_{ix})^2W^\varepsilon_{i\alpha \alpha}=0,
\end{equation}
if we choose $\alpha^\varepsilon_i$ to satisfy
\begin{equation}
\alpha^\varepsilon_{it}+\tilde{\lambda_i}\alpha^\varepsilon_{ix}-\varepsilon\alpha^\varepsilon_{ixx}=0.
\end{equation}
We impose the following initial conditions
\begin{equation}
t=0, \alpha^\varepsilon_i=x, W_i^\varepsilon=R_i(u_{0\varepsilon}(x)).
\end{equation}
Let $\Theta^\varepsilon_i=\alpha^\varepsilon_{ix}$, then,
\begin{equation}
\Theta^\varepsilon_{it}+(\tilde{\lambda_i}\Theta^\varepsilon_{i})_x-\varepsilon\Theta^\varepsilon_{ixx}=0,
\end{equation}
\begin{equation}
\Theta^\varepsilon _i (0,x)=1.
\end{equation}
By Maximum principle, $\Theta^\varepsilon_i>0$,which implies $\alpha^\varepsilon (t,x)$ is an increasing function of x. It also follows from \eqref{4.10} that
\begin{equation}\label{4.19}
|\partial_\alpha W^\varepsilon_i(t,\alpha)|_{L^p(R)}\le |\partial_xR_i(u_{0\varepsilon}(x))|_{L^p(R)}\lesssim M.
\end{equation}

We start the proof of Theorem 1.1 with the parabolic estimate. By \eqref{4.1},\eqref{4.2}, we get
\begin{equation}\label{4.20}
u(t,x)=E(t,x)\ast J_{\varepsilon}\ast u_{0}-\int_0^T E(t-\tau, x)\ast (A(u)u_x)(\tau) d\tau
\end{equation}
where
\begin{equation}\label{4.21}
E(t,x)=\frac{1}{\sqrt{2\pi\varepsilon t}}\exp{-\frac{x^2}{4\varepsilon t}}
\end{equation}
is the heat kernel. We take $\delta=\delta (D)$ to be a constant depending only on D and $\delta\ll D$.
We consider first the equation on the time interval $t\in [0, \varepsilon\delta]$, by \eqref{4.20}, we have
\begin{equation}
u_x(t,x)=\varepsilon^{-1}E(t,x)\ast J'_{\varepsilon}\ast u_{0}-\int_0^T E_x(t-\tau, x)\ast (A(u)u_x)(\tau) d\tau.
\end{equation}
Therefore, we get
\begin{equation}
|u_x(t)|_{L^\infty (R)}\lesssim \varepsilon^{-1} |u_{0}|_{L^\infty(R)}+\int_0^T \frac{1}{\sqrt{\varepsilon (t-\tau)}}|u_x(\tau)|_{L^\infty (R)} d\tau.
\end{equation}
Thus, we have
\begin{equation}
\sup_{0\le t\le \varepsilon\delta }|u_x(t)|_{L^\infty(R)}\lesssim  \varepsilon^{-1}+\sqrt{\delta}\sup_{0\le t\le \varepsilon\delta }|u_x(t)|_{L^\infty(R)},
\end{equation}
which implies
\begin{equation}\label{4.24}
\sup_{0\le t\le \varepsilon\delta }|u_x(t)|_{L^\infty(R)}\lesssim  \varepsilon^{-1},
\end{equation}
provided that $\delta$ is taken to be small enough.
Now, let $t\ge \varepsilon\delta$, we consider the equation on the time interval $s\in [t-\varepsilon\delta, t]$, we have
\begin{equation}\label{4.25}
u(s,x)=E(s-(t-\varepsilon\delta),x)\ast u(t-\varepsilon)-\int_{t-\varepsilon\delta }^s E(s-\tau, x)\ast (A(u)u_x)(\tau) d\tau,
\end{equation}
and thus,
\begin{equation}\label{4.26}
u_x(s,x)=E_x(s-(t-\varepsilon\delta),x)\ast u(t-\varepsilon\delta)-\int_{t-\varepsilon\delta }^s E_x(s-\tau, x)\ast (A(u)u_x)(\tau) d\tau.
\end{equation}
Therefore, we get
\begin{equation}
|u_x(s)|_{L^\infty (R)}\lesssim \frac{1}{\sqrt{\varepsilon (s-t+\varepsilon\delta)}} |u(t-\varepsilon\delta )|_{L^{\infty}}+\int_{t-\varepsilon\delta }^s \frac{1}{\sqrt{\varepsilon (s-\tau)}}|u_x(\tau)|_{L^\infty (R)} d\tau.
\end{equation}
Thus, we have
\begin{equation}
\sqrt{ (s-t+\varepsilon\delta)}|u_x(s)|_{L^\infty (R)}\lesssim \varepsilon^{-1/2}+\varepsilon^{-1/2}\sqrt{ (s-t+\varepsilon\delta)}\sup_{t-\varepsilon\delta\le s\le t}
(\sqrt{ (s-t+\varepsilon\delta)}|u_x(s)|_{L^\infty (R)})
\end{equation}
which implies
\begin{equation}
\sup_{t-\varepsilon\delta\le s\le t}(\sqrt{ (s-t+\varepsilon\delta)}|u_x(s)|_{L^\infty (R)})\lesssim \varepsilon^{-1/2}
\end{equation}
provided that $\delta $ is taken to be sufficiently small. Take $s=t$, we get
\begin{equation}\label{4.30}
|u_x(t)|_{L^\infty (R)}\lesssim \varepsilon^{-1}\quad \forall t\ge\varepsilon\delta .
\end{equation}
By \eqref{4.24} and \eqref{4.30}, we finally arrive at
\begin{equation}\label{4.31}
|u_x(t)|_{L^\infty (R)}\lesssim \varepsilon^{-1}\quad \forall t\ge 0.
\end{equation}
It follows then $|\tilde{\lambda}_i|\le C_3$ for some constant $C_3$ depending only on $D$. Thus, by maximum principle
\begin{equation}
x-C_3t\le \alpha^\varepsilon_i(t,x)\le x+C_3t
\end{equation}
because the left hand side is a subsolution and the right hand side is a supsolution. Thus, there exists a subsequence such that $\alpha_i^\varepsilon (t,x)$ converges to some $\alpha_i(t,x)$ all most everywhere. By \eqref{4.19}, there exists a further subsequence such that $W_i^\varepsilon$ weakly converges to some $W_i(t,\alpha)$ in $\dot{W}^{1,p}$, by Sobolev embedding Theorem, $W_i^\varepsilon$ strongly converges to $W_i(t,\alpha)$ in Holder space. Therefore, taking limit in
\begin{equation}
u^\varepsilon(t,x)=U(W_1^\varepsilon (t,\alpha^\varepsilon_1(t,x)),W_2^\varepsilon (t,\alpha^\varepsilon_2(t,x)) ),
\end{equation}
we conclude the proof of Theorem 1.1.

\section*{Acknowledgement}

The author is supported by Key Laboratory of Mathematics for Nonlinear Sciences (Fudan University), Ministry of Education of China, P.R.China. Shanghai Key Laboratory for Contemporary Applied Mathematics, School of Mathematical Sciences, Fudan University, P.R. China, and by Shanghai Science and Technology Program [Project No. 21JC1400600].

\bibliographystyle{plain}
%\bibliography{bibli}
\end{CJK*}

\end{document}